\documentclass[11pt]{amsart}
\usepackage[usenames,dvipsnames]{pstricks}
\usepackage[mathscr]{eucal}
\usepackage{amsfonts,amsmath,amssymb,amsthm,amscd,amsxtra}
\usepackage{enumerate,verbatim}
\usepackage[all,2cell,ps]{xy}
\usepackage[notcite, notref]{}
\usepackage[pagebackref]{hyperref}
\usepackage{mathptmx}
\usepackage{calc,color}
\usepackage{wasysym}

\theoremstyle{plain} 

\newtheorem{thm}{Theorem}[section]
\newtheorem{dfn}[thm]{Definition}

\newtheorem{prop}[thm]{Proposition}

\newtheorem{cor}[thm]{Corollary}
\newtheorem{lem}[thm]{Lemma}

\numberwithin{equation}{section}

\newcommand{\fm}{\mathfrak{m}}

\newcommand{\fp}{\mathfrak{p}}
\newcommand{\fq}{\mathfrak{q}}
\newcommand{\fa}{\mathfrak{a}}
\newcommand{\fb}{\mathfrak{b}}
\newcommand{\fc}{\mathfrak{c}}

\newtheorem{chunk}[thm]{\hspace*{-1.065ex}\bf}

\DeclareMathOperator{\ann}{ann} \DeclareMathOperator{\Ass}{Ass}
 
\DeclareMathOperator{\V}{V} \DeclareMathOperator{\hh}{H}
\DeclareMathOperator{\E}{E} \DeclareMathOperator{\Att}{Att}
 \DeclareMathOperator{\X}{X}
\DeclareMathOperator{\gor}{Gor}
\DeclareMathOperator{\cc}{c}
\def\gd{\operatorname{\mathsf{G-dim}}}

\def\gkd{\operatorname{\mathsf{G}_{\it C}\mathsf{-dim}}}
\def\gfcd{\operatorname{\mathsf{G}_{\it \fc}\mathsf{-dim}}}
\def\gfcdp{\operatorname{\mathsf{G}_{\it {\fc}_\fp}\mathsf{-dim}}}
\def\gkkd{\operatorname{\mathsf{G}_{\it K}\mathsf{-dim}}}
\def\gkkkd{\operatorname{\mathsf{G}_{\it K}}}
\def\gc{\operatorname{\mathsf{G}_{\it C}}}
\def\gfc{\operatorname{\mathsf{G}_{\it \fc}}}

\def\pd{\operatorname{\mathsf{pd}}}

\def\gr{\operatorname{\mathsf{grade}}}

\def\Tr{\mathsf{Tr}}
\def\trk{\mathsf{Tr}_{C}}

\def\Min{\mathsf{Min}}

\DeclareMathOperator{\coker}{Coker}

\def\depth{\operatorname{\mathsf{depth}}}
 
\def\Ext{\operatorname{\mathsf{Ext}}}

\def\Hom{\operatorname{\mathsf{Hom}}}

\DeclareMathOperator{\id}{id} 

\DeclareMathOperator{\Supp}{Supp} \DeclareMathOperator{\Spec}{Spec}

\def\Tor{\operatorname{\mathsf{Tor}}}

\def\urltilda{\kern -.15em\lower .7ex\hbox{\~{}}\kern .04em}
\def\urldot{\kern -.10em.\kern -.10em}\def\urlhttp{http\kern -.10em\lower -.1ex
\hbox{:}\kern -.12em\lower 0ex\hbox{/}\kern -.18em\lower
0ex\hbox{/}}

\begin{document}
\baselineskip=15pt

\title[Notes on linkage of modules]
 {Notes on linkage of modules}

\bibliographystyle{amsplain}
\author[A. Sadeghi]{Arash Sadeghi}

\address{School of Mathematics, Institute for Research in Fundamental Sciences (IPM), P.O. Box: 19395-5746, Tehran, Iran }
\email{sadeghiarash61@gmail.com}

\keywords{linkage of modules, attached prime, local cohomology\\
This research was supported by a grant from IPM}

 \maketitle
\begin{abstract}
	Let $R$ be a Cohen-Macaulay local ring. It is shown that under some mild conditions, the Cohen-Macaulayness property is preserved under linkage. We also study the connection of $(S_n)$ locus of a horizontally linked module and the attached primes of certain local cohomology modules of its linked module.
    
\end{abstract}

\section{introduction}
The theory of linkage for subschemes of projective space goes back more than a century in some sense, but the modern study was introduced by Peskine and Szpiro \cite{PS} in 1974. Recall that two ideals $I$ and $J$ in a Cohen-Macaulay local ring $R$ are said to be linked if there is a regular sequence $\alpha$ in their intersection such that $I=(\alpha:J)$ and $J=(\alpha:I)$. The first main result in the theory of linkage, due to Peskine and Szpiro, indicates that the Cohen-Macaulayness property is preserved under linkage over Gorenstein rings. They also give a counterexample to show that the above result is no longer true if the base ring is Cohen-Macaulay but not Gorenstein. On the other hand, there are interesting extensions of the Peskine-Szpiro Theorem
to the case that $R$ is Cohen-Macaulay and one of the ideals $I$ or $J$ is strongly Cohen-Macaulay (i.e. Koszul homology modules are Cohen-Macaulay) \cite{Hu} or satisfies
the sliding depth condition (on Koszul homology) \cite{HVV}. In a different direction, Schenzel generalized the Peskine-Szpiro Theorem to the vanishing of certain local cohomology modules. More precisely, for linked ideals $I$ and $J$ over a Gorenstein local ring $R$ with maximal ideal $\fm$, the Serre condition $(S_n)$ for $R/I$ is equivalent to the vanishing of
the local cohomology modules $\hh^i_{\fm}(R/J)= 0$ for all $i$, $\dim(R/J)-n<i<\dim(R/J)$ \cite[Theorem 4.1]{Sc}.

A few authors advanced the theory to the setting of linkage of modules in different ways, for instance  Martin \cite{Mar}, Yoshino and Isogawa \cite{YI}, Martsinkovsky and
Strooker \cite{MS}, and Nagel \cite{Na}. Based on these generalizations, several works have been done on studying the linkage theory in the context of modules; see for example \cite{CDGST}--\cite{DS2}, \cite{IT}, \cite{N}, \cite{TJP} and \cite{TP}.
In this paper, we are interested in linkage of modules in the sense of \cite{MS}. Martsinkovsky and Strooker generalized the notion of linkage for modules by using the composition of two functors: transpose and syzygy. They  showed that ideals $\fa$ and $\fb$ are linked by zero ideal if and only if $R/\fa\cong\lambda(R/\fb)$ and $R/\fb\cong\lambda(R/\fa)$, where $\lambda:=\Omega\Tr$.

Recall that the $(S_k)$ locus of $M$, denoted by $S_k(M)$, is  the subset of $\Spec R$ consisting of all prime ideals $\fp$ of $R$ such that $M_\fp$ satisfies the Serre condition $(S_k)$. The Gorenstein locus of $R$, denoted by $\gor(R)$, is the subset of $\Spec R$ consisting of all prime ideals $\fp$ of $R$ such that $R_\fp$ is a Gorenstein local ring.
In the fisrt part of this paper, we study the connection of $(S_n)$ locus of a horizontally linked module and the attached primes of certain local cohomology modules of its linked module. As a consequence, we obtain the following result (see also Theorem \ref{t3} for a more general case).\\

\textbf{Theorem A.}
{\it Let $(R,\fm)$ be a Cohen-Macaulay local ring of dimension $d$ and let $X$ be an open subset of $\Spec R$ which is contained in $\gor(R)$. Assume that $M$ is a horizontally linked $R$--module and that $n$ is a positive integer. Then $X\subseteq S_n(M)$ if and only if $\Att_R(\hh^i_{\fm}(\lambda M))\subseteq\Spec R\setminus X$ for $d-n<i<d$.} \\

Note that for an Artinian $R$--module $N$, $\Att_R(N)=\emptyset$ if and only if $N=0$. Therefore, the above theorem can be viewed as a generalization of Schenzel's result. We should note that Theorem A is new even in the case of cyclic module. The above theorem helps us to study the attached primes of local cohomology modules of a linked module. As an application of Theorem A, for a horizontally linked module $M$ over a Gorenstein local ring $R$ and a closed subset $X$ of $\Spec R$, it is shown that
$\Att_R(\hh^i_\fm( M))\subseteq X$ for $0<i<\dim R$ if and only if $\Att_R(\hh^i_\fm(\lambda M))\subseteq X$ for $0<i<\dim R$ (see Corollary \ref{c}). 

There is another interesting generalization of the Peskine-Szpiro Theorem, due to Martin \cite{Mar1}. Let $R$ be a Gorenstein local ring and let $I$ be a Cohen-Macaulay ideal linked to $J$ by an ideal $\fa$, where $\fa$ is Cohen-Macaulay and generically Gorenstein. Then $J$ is Cohen-Macaulay if and only if $I+\fb$ is a Cohen-Macaulay ideal with the same height as $\fb$, where $\fb$ is an ideal of $R$, containing $\fa$, with $\fb/\fa$ equal to a canonical module of $\omega_{R/\fa}$ \cite[Corolary 2.11]{Mar1}.
In the second part of this paper, we study the theory of linkage for Cohen-Macaulay modules. As a consequence, 
we generalize the Martin's result to the module case. 
More precisely, we have the following result (see also Corollary \ref{ttttt} for a more general case). \\

\textbf{Theorem B.}
{\it	Let $R$ be a Cohen-Macaulay local ring with canonical module $\omega_R$ and let $\fa$ be a generically Gorenstein, Cohen-Macaulay ideal of $R$. Write $\omega_{R/\fa}\cong\fc/\fa$. Assume that $M$ is a Cohen-Macaulay $R$--module which is linked to an $R$--module $N$ by the ideal $\fa$. Then $N$ is Cohen-Macaulay if and only if $\depth_R(M/\fc M)\geq\dim(R/\fa)-1$.
}\\

The organization of the paper is as follows. In Section 2, we collect preliminary notions, definitions and some known results which will be used in this paper. In Section 3, we 
study the connection of $(S_n)$ locus of a horizontally linked module and the attached primes of certain local cohomology modules of its linked module. As a consequence,  we  prove Theorem A in this section. In section 4, we study the theory of linkage for Cohen-Macaulay modules. We prove Theorem B in this section. In section 5, we establish a duality on local cohomology modules of a linked module. 


\section{Definitions and preliminary results}
Throughout the paper $R$ is a commutative Noetherian local ring and all modules over $R$ are assumed to be
finitely generated. We start by recalling
several definitions and terminology from \cite{AB,BS,MS}.

\begin{chunk}[\cite{AB}] \label{AuBrsequence} 
	Let $M$ be an $R$--module. For a positive integer $i$, we denote by $\Omega^iM$ the $i$-th syzygy of M, namely, the image of the $i$-th differential map in a minimal free resolution of $M$.
	
	The \emph{transpose} $\Tr M$ of $M$ is defined as the cokernel of the $R$-dual map $\partial_1^{\ast}=\Hom_{R}(\partial_1,R)$ of the first differential map $\partial_1$ in a minimal free resolution of $M$. Hence there is an exact sequence of the form
	$0\rightarrow M^*\rightarrow P_0^*\rightarrow P_1^*\rightarrow \Tr M\rightarrow 0$. Note that the modules $\Omega^iM$ and $\Tr M$ are uniquely determined up to isomorphism, since so is a minimal free resolution of $M$, and there is a stable isomorphism $\Omega^2\Tr M \cong M^{\ast}$. 
	
	A \emph{stable} module is a module with no non-zero
	free direct summands. An $R$--module $M$ is called a \emph{syzygy module} if it is
	embedded in a free $R$--module.
\end{chunk}
\begin{chunk}\label{link}\textbf{Linkage of modules.}(\cite{MS})
	Two $R$--modules $M$ and $N$ are said
	to be \emph{horizontally linked} if $M\cong \lambda N$ and
	$N\cong\lambda M$, where $\lambda M:=\Omega\Tr M$. Equivalently, $M$ is
	horizontally linked (to $\lambda M$) if and only if $M\cong\lambda^2 M$. An $R$--module $M$ is horizontally linked if and only if it is stable and $\Ext^1_R(\Tr M,R)=0$, equivalently $M$ is stable and is a syzygy module \cite[Theorem 2 and Proposition 3]{MS}.
	An $R$--module
	$M$ is said to be \emph{linked} to an $R$--module $N$, by an ideal $\fc$ of $R$, if $\fc\subseteq\ann_R(M)\cap\ann_R(N)$ and $M$ and
	$N$ are horizontally
	linked as $R/{\fc}$--modules. In this situation we write $M\underset{\fc}{\sim}N$.
	
	Let $M$ be a stable $R$--module and let $P_1\rightarrow P_0\rightarrow M\rightarrow 0$ be a minimal projective presentation of $M$. Then
	$P_0^*\rightarrow P_1^*\rightarrow \Tr M\rightarrow 0$ is a minimal projective presentation of $\Tr M$ \cite[Theorem 32.13]{AF}. We quote the following induced exact sequences
	\begin{equation}\tag{\ref{link}.1}
	0\longrightarrow M^*\longrightarrow P_0^*\longrightarrow \lambda M\longrightarrow0,
	\end{equation}
	\begin{equation}\tag{\ref{link}.2}
	0\longrightarrow \lambda M\longrightarrow P_1^*\longrightarrow \Tr M\longrightarrow0.
	\end{equation}
\end{chunk}
\begin{chunk} \textbf{Gorenstein dimension with respect to a semidualizing module.} \label{Gdim} (\cite{F}, \cite{G})
	An $R$--module $C$ is called a semidualizing module, if the homothety morphism $R\to\Hom_R(C,C)$ is an isomorphism and $\Ext^i_R(C,C)=0$ for all $i>0$. Semidualizing modules are initially studied by Foxby \cite{F} and Golod \cite{G}. It is clear that $R$ itself is a semidualizing $R$--module. Over a Cohen-Macaulay local ring $R$, a canonical module $\omega_R$ of $R$ is a semidualizing module with finite injective dimension. 
	
	An $R$--module $M$ is said to be \emph{totally $C$-reflexive} if the
	natural map $M \to \Hom(\Hom(M,C),C)$ is bijective and $\Ext_{R}^{i}(M,C)=0=\Ext_{R}^{i}(\Hom(M,C),C)$ for all $i\geq 1$. The smallest nonnegative integer $n$ for which there exists an exact sequence
	$$0 \to X_{n} \to X_{n-1} \to \dots \to X_{1} \to X_{0} \to M\to 0,$$
	such that each $X_{i}$ is totally $C$-reflexive, is called the \emph{$\gc$-dimension} of $M$. If $M$ has $\gc$-dimension $n$, we write $\gkd_R(M)=n$. Therefore $M$ is totally $C$-reflexive if and only if $\gkd_R(M)=0$.
	
	Let $P_1\overset{f}{\rightarrow}P_0\rightarrow M\rightarrow 0$ be a
	projective presentation of $M$. The transpose of $M$ with respect to
	$C$, $\trk M$, is defined to be $\coker f^{\triangledown}$, where
	$(-)^{\triangledown} := \Hom_R(-,C)$, which satisfies the exact
	sequence
	\begin{equation}\tag{\ref{Gdim}.1}
	0\rightarrow M^{\triangledown}\rightarrow
	P_0^{\triangledown}\rightarrow P_1^{\triangledown}\rightarrow \trk
	M\rightarrow 0.
	\end{equation}
	
	We catalogue some properties of Gorenstein dimension with respect to a semidualizing module:
	
	\begin{enumerate}[\rm(i)]
		\item If $C$ is a dualizing $R$--module, then $\gkd_R(M)<\infty$ for all $R$--modules $M$;
		\item If $\gkd_R(M)<\infty$, then $\gkd_R(M)=\depth(R)-\depth(M)$;
		\item If $\gkd_R(M)<\infty$, then $\gkd_R(M)=\sup\{i\mid\Ext^i_R(M,C)\neq0\};$
		\item $\gkd_R(M)=0$ if and only if $\gkd_R(\trk M)=0$.
	\end{enumerate}
		An $R$--module $M$ is said to satisfy the property $\widetilde{S}_k$ if $\depth_{R_\fp}(M_{\fp})\geq\min\{k, \depth R_\fp\}$ for all $\fp\in\Spec R$. Recall that the $\widetilde{S}_k$ locus of $M$, denoted by $\widetilde{S}_k(M)$, is defined by
	$$\widetilde{S}_k(M)=\{\fp\in\Spec R\mid M_\fp \text{ satisfies } \widetilde{S}_k\}.$$
	
	For a horizontally linked module $M$ over a
	Cohen-Macaulay local ring $R$, the properties $\widetilde{S}_k$ and
	$(S_k)$ are identical. An $R$--module $M$ is called $n$-torsionfree, when $\Ext^i_R(\Tr M,R)=0$ for $1\leq i\leq n$.
	
	The following Theorem is due to Auslander and Bridger.
	\begin{thm}\cite[Theorem 4.25]{AB}\label{AS}
		Let $M$ be an $R$--module of finite Gorenstein dimension. Then $M$ is $n$-torsionfree if and only if $M$ satisfies $\widetilde{S}_n$.	
	\end{thm}
\end{chunk}
\begin{chunk}\label{ac}\textbf{Auslander class.}
		Let $C$ be a semidualizing $R$-module. The Auslander class with respect to $C$, denoted by $\mathcal{A}_C$, consists of all $R$--modules $M$ satisfying the following conditions.
	\begin{enumerate}[(i)]
		\item{The natural map $\mu:M\to\Hom_R(C, M\otimes_RC)$ is an isomorphism.}
		\item{$\Tor^R_i(M,C) =0=\Ext^i_R(C, M\otimes_RC)$ for all $i>0$.}
	\end{enumerate}
	Every module of finite projective dimension is in the Auslander class $\mathcal{A}_C$. Over a Cohen-Macaulay local ring $R$ with canonical module $\omega_R$, we have $M\in\mathcal{A}_{\omega_R}$ if and only if $\gd_R(M)<\infty$ \cite[Theorem 1]{F1}.
\end{chunk}
\begin{chunk}\textbf{Secondary representation.}(\cite{IGM}) Let $M$ be a nonzero Artinian $R$--module. We say that $M$ is secondary if the multiplication by $x$ on $M$ is surjective or nilpotent for every $x\in R$. In this case, $\fp:=\sqrt{\ann_R(M)}$ is a
prime ideal of $R$ and we say that $M$ is $\fp$-secondary. Note that every Artinian $R$--module $M$ has a minimal secondary representation $M=M_1+\cdots+M_n$, where $M_i$ is $\fp_i$-secondary, each $M_i$ is not redundant and $\fp_i\neq\fp_j$ for all $i\neq j $. The set $\{\fp_1,\cdots,\fp_n\}$ is independent of the choice of the minimal secondary representation of $M$. This set is called the set of attached primes of $M$ and denoted by $\Att_R(M)$. 

We denote by $\widehat{M}$ the completion of $M$ in the $\fm$-adic topology.
In the following we collect some basic properties of attached primes of a module which will be used throughout the paper (see \cite{IGM} and \cite{BS}).
\begin{thm}\label{att}
Let $M$ be an Artinian $R$--module. The following statements hold.
\begin{enumerate}[(i)]
	\item{$M\neq0$ if and if $\Att_R(M)\neq\emptyset.$}
	\item{$M$ has finite length if and only if $\Att_R(M)\subseteq\{\fm\}$.}
	\item{$\Att_R(M)=\Ass_R(\Hom_R(M,\E_R(k)))$, where $\E_R(k)$ is the injective envelope of $k$.}
	\item{$\Att_R(M)=\{P\cap R\mid P\in\Att_{\widehat{R}}(M)\}$.}
\end{enumerate}	
\end{thm}
The following result plays an important role in the proof of Theorem A.
\begin{thm}\cite[Theorem 1.1]{NQ}\label{nqat}
	Let $R$ be a ring which is the homomorphic image of a Cohen-Macaulay local ring. Then
	$$\Att_{R_\fp}(\hh^{i-\dim R/\fp}_{\fp R_\fp}(M_\fp))=\{\fq R_\fp\mid \fq\in\Att_{R}(\hh^i_{\fm}(M), \fq\subseteq\fp)\},$$
	for every $R$--module $M$, integer $i\geq0$ and prime ideal $\fp\in R$.
\end{thm}
\end{chunk}


\section{Linkage and the attached primes of local cohomology modules}
Recall that a subset $X$ of $\Spec R$ is called \emph{stable under generalization} provided that if $\fp\in X$ and
$\fq\in\Spec R$ with $\fq\subseteq\fp$, then $\fq\in X$. Note that every open subset of $\Spec R$ is stable under generalization. For an integer $n$, set $X^n(R):=\{\fp\in\Spec(R)\mid \depth R_\fp\leq n\}$.

We start by recalling the following result.
\begin{lem}\label{l1}
Let $C$ be a semidualizing $R$--module and let $M$ be an $R$--module. 
Assume that $M\in\mathcal{A}_C$ and that $n$ is a positive integer. Consider the following statements.
\begin{enumerate}[(i)]
\item{$\Ext^i_R(\Tr M,R)=0$ for $1\leq i\leq n$.}
\item{$\Ext^i_R(\Tr M,C)=0$ for $1\leq i\leq n$.}
\item{$M$ satisfies $\widetilde{S}_n$.}
\end{enumerate}
Then we have the following.
\begin{enumerate}[(a)]
	\item{(i)$\Leftrightarrow$(ii)$\Rightarrow$(iii)}
	\item{If $\gd_{R_\fp}(M_\fp)<\infty$ for all
		$\fp\in\X^{n-1}(R)$ (e.g. $\id_{R_\fp}(C_\fp)<\infty$ for all	$\fp\in\X^{n-1}(R)$), then all the statements (i)-(iii) are equivalent.} 
\end{enumerate}
\end{lem}
\begin{proof}
See \cite[Theorem 4.1]{Sa} and \cite[Theorem 2.12]{DS1}.	
\end{proof}
The proof of Theorem A is based on the following lemma, which is of independent interest.
\begin{lem}\label{l2}
	Let $C$ be a semidualizing $R$--module and let $M$ be an $n$-torsionfree $R$--module.
	Assume that $M\in\mathcal{A}_C$ and that $\gd_{R_\fp}(M_\fp)<\infty$ for all $\fp\in X^{n}(R)$ (e.g. $\id_{R_\fp}(C_\fp)<\infty$ for all $\fp\in X^{n}(R)$). Then
	$\Ass_R(\Ext^{n+1}_R(\Tr M,C))=\Ass_R(\Ext^{n+1}_R(\Tr M,R)).$
\end{lem}
\begin{proof}
	We argue by induction on $n$. Let $n=0$ and let $\fp\in\Ass(\Ext^{1}_R(\Tr M,R))$. Consider the following exact sequence (see \cite[Proposition 5]{M1})
	\begin{equation}\tag{\ref{l2}.1}
	0\rightarrow\Ext^1_R(\Tr M,R)\rightarrow M\rightarrow M^{**}\rightarrow\Ext^2_R(\Tr M,R)\rightarrow0.
	\end{equation}
	It follows from the exact sequence (\ref{l2}.1) that $\fp\in\Ass_R(M)$. As $M\in\mathcal{A}_C$, we have $M\cong\Hom_R(C,C\otimes_RM)$. Therefore $\fp\in\Ass_R(M\otimes_RC)$.
	Let us now consider the following exact sequence (see \cite[Proposition 2.6]{AB})
	\begin{equation}\tag{\ref{l2}.2}
	0\rightarrow\Ext^1_R(\Tr M,C)\rightarrow M\otimes_RC\rightarrow\Hom_R(M^*,C)\rightarrow\Ext^2_R(\Tr M,C)\rightarrow0.
	\end{equation}
	If $\depth R_\fp=0$ then $\gd_{R_\fp}(M_\fp)=0$ and so $\Ext^1_R(\Tr M,R)_\fp=0$ which is a contradiction.
	Therefore, $\depth R_\fp>0$ and so $\depth_{R_\fp}(\Hom_R(M^*,C)_\fp)\geq\min\{2,\depth R_\fp\}>0$. Hence, it follows from the exact sequence
	(\ref{l2}.2) that $\depth_{R_\fp}(\Ext^1_R(\Tr M,C)_\fp)=0$. In other words, $\fp\in\Ass_R(\Ext^1_R(\Tr M,C))$. Similarly, one can prove that
	$\Ass_R(\Ext^{1}_R(\Tr M,C))\subseteq\Ass_R(\Ext^{1}_R(\Tr M,R))$. Now let $n>0$. As $M$ is $n$-torsionfree, $\Ext^i_R(\Tr M,R)=0=\Ext^i_R(\Tr M,C)$ for $1\leq i\leq n$ by Lemma \ref{l1}.
	Consider the universal pushforward of $M$:
	\begin{equation}\tag{\ref{l2}.3}
	0\longrightarrow M\longrightarrow F\longrightarrow M_1\longrightarrow0,
	\end{equation}
	where $F$ is free and $\Ext^1_R(M_1,R)=0$ (see \cite{M1}). Hence, by \cite[Lemma 3.9]{AB}, the exact sequence (\ref{l2}.3) induces the following exact sequence
	\begin{equation}\tag{\ref{l2}.4}
	0\longrightarrow \Tr M_1\longrightarrow P\longrightarrow \Tr M\longrightarrow0,
	\end{equation}
	where $P$ is a free $R$-module. It follows from the exact sequence (\ref{l2}.4) that $M_1$ is $(n-1)$-torsionfree.
	Also, by the exact sequence (\ref{l2}.3), $M_1\in\mathcal{A}_C$.
	Hence, by induction hypothesis
	\[\begin{array}{rl}
	\Ass_R(\Ext^{n+1}_R(\Tr M,R))& =\Ass_R(\Ext^{n}_R(\Tr M_1,R))\\
	&=\Ass_R(\Ext^n_R(\Tr M_1,C))=\Ass_R(\Ext^{n+1}_R(\Tr M,C)).
	\end{array}\]
\end{proof}
Now we can present our first main result.
\begin{thm}\label{t3}
Let $(R,\fm)$ be a formally equidimensional local ring of dimension $d$ which is the homomorphic image of a Cohen-Macaulay ring. Assume that $M$ is a horizontally linked $R$--module and that $n$ is a positive integer. Assume further that $X$ is a subset of $\Spec R$ which is stable under generalization and that $R_{\fp}$ is Cohen-Macaulay and $\gd_{R_\fp}(M_\fp)<\infty$ for all $\fp\in X$ (e.g. $X\subseteq\gor(R)$).
Then the following are equivalent.
\begin{enumerate}[(i)]
\item{$X\subseteq\widetilde{S}_n(M)$.}
\item{$\Att_R(\hh^i_{\fm}(\lambda M))\subseteq\Spec R\setminus X$ for $d-n<i<d$.}
\end{enumerate}
\end{thm}
\begin{proof}
(i)$\Rightarrow$(ii). Assume contrarily that  $\fp\in\Att_R(\hh^i_\fm(\lambda M))\cap X$ for some $d-n<i<d$. Therefore,
$\fp R_\fp\in\Att_{R_\fp}(\hh^{i-\dim R/\fp}_{\fp R_\fp}(\lambda_{R_\fp} M_\fp))$ by Theorem \ref{nqat}.
As $\fp\in\widetilde{S}_n(M)\cap X$, we have $\hh^{j}_{\fp R_\fp}(\lambda_{R_\fp} M_\fp)=0$ for $\dim R_\fp-n<j<\dim R_\fp$ by \cite[Theorem 4.2]{DS}.
Note that by \cite[Theorem 31.5]{Ma}, $R$ is equidimensional and catenary. Hence
$d=\dim R_\fp+\dim R/\fp$ and so $\dim R_\fp-n<i-\dim R/\fp<\dim R_\fp$.
In particular, $\hh^{i-\dim R/\fp}_{\fp R_\fp}(\lambda_{R_\fp} M_\fp)=0$ which is a contradiction by Theorem \ref{att}(i).

(ii)$\Rightarrow$(i). Assume contrarily that $X\nsubseteq\widetilde{S}_n(M)$. Set
$Y=\{\fp\in X\mid M_\fp \text{ does not satisfy } \widetilde{S}_n \}$. Let $\fp_0\in\Min(Y)$. As $\gd_{R_{\fp_0}}(M_{\fp_0})<\infty$, $\Ext^i_R(\Tr M,R)_{\fp_0}\neq0$
for some $1\leq i\leq n$ by Theorem \ref{AS}. Set
\begin{equation}\tag{\ref{t3}.1}
t=\min\{j\mid\Ext^j_R(\Tr M,R)_{\fp_0}\neq0, 1\leq j\leq n\}.
\end{equation}
As $M$ is horizontally linked, $\Ext^1_R(\Tr M,R)=0$ and so $1<t\leq n$.
Note that $\fq\in\widetilde{S}_n(M)$ for all prime ideal $\fq$ with $\fq\subsetneq\fp_0$, because $\fp_0\in\Min(Y)$ and $X$ is stable under generalization.
In other words, by Theorem \ref{AS},
\begin{equation}\tag{\ref{t3}.2}
\Ext^t_R(\Tr M,R)_\fq=0 \text{ for all prime ideal } \fq \text{ with } \fq\subsetneq\fp_0.
\end{equation}
Therefore $\fp_0$ is a minimal element of $\Supp_R(\Ext^t_R(\Tr M,R))$ and so $\fp_0\in\Min_R(\Ext^t_R(\Tr M,R))$.
Note that by \cite[Proposition 9.A]{Ma1} $$\Ass_R(\Ext^t_R(\Tr M,R))=\{P\cap R\mid P\in\Ass_{\widehat{R}}(\Ext^t_{\widehat{R}}(\widehat{\Tr M},\widehat{R}))\}.$$
Let $P_0\in\Min_{\widehat{R}}(\Ext^t_{\widehat{R}}(\widehat{\Tr M},\widehat{R}))$ such that $P_0\cap R=\fp_0$.
It follows easily from (\ref{t3}.1) that
\begin{equation}\tag{\ref{t3}.3}
\Ext^i_{\widehat{R}}(\widehat{\Tr M},\widehat{R})_{P_0}=0 \text{ for all } 0<i<t.
\end{equation}
As $\gd_{R_{\fp_0}}(M_{\fp_0})<\infty$ and the local homomorphism $R_{\fp_0}\rightarrow \widehat{R}_{P_0}$ is flat,
$\gd_{\widehat{R}_{P_0}}(\widehat{M}_{P_0})<\infty$. As $R$ is the homomorphic image of a Cohen-Macaulay ring, all formal fibre of $R$ is Cohen-Macaulay.
Therefore, $\widehat{R}_{P_0}/{\fp_0\widehat{R}_{P_0}}\cong (R_{\fp_0}/{\fp_0 R_{\fp_0}}\otimes_R\widehat{R})_{P_0}$ is Cohen-Macaulay and so $\widehat{R}_{P_0}$ is Cohen-Macaulay (see \cite[Page 181]{Ma}).
Therefore, by Lemma \ref{l2}, (\ref{t3}.3) and (\ref{ac}),
\begin{equation}\tag{\ref{t3}.5}
P_0\widehat{R}_{P_0}\in\Ass(\Ext^t_{\widehat{R}_{P_0}}(\Tr_{\widehat{R}_{P_0}}\widehat{M}_{P_0},\omega_{\widehat{R}_{P_0}})).
\end{equation}
As $t>1$, $\Ext^t_{\widehat{R}_{P_0}}(\Tr_{\widehat{R}_{P_0}}\widehat{M}_{P_0},\omega_{\widehat{R}_{P_0}})\cong\Ext^{t-1}_{\widehat{R}_{P_0}}(\widehat{\lambda M}_{P_0},\omega_{\widehat{R}_{P_0}})$.
Therefore, by (\ref{t3}.5), Theorem \ref{att}(iii) and the local duality theorem,  $P_0\widehat{R}_{P_0}\in\Att_{\widehat{R}_{P_0}}(\hh^{s-t+1}_{P_0\widehat{R}_{P_0}}(\widehat{\lambda M}_{P_0}))$ where $s=\dim \widehat{R}_{P_0}$. As $\widehat{R}$ is equidimensional, $s+\dim(\widehat{R}/{P_0})=d$. Hence, by Theorem \ref{nqat}, $P_0\in\Att_{\widehat{R}}(\hh^{d-t+1}_{\widehat{\fm}}(\widehat{\lambda M}))$ and so
$\fp_0\in\Att_R(\hh^{d-t+1}_{m}(\lambda M))$ by Theorem \ref{att}(iv), which is a contradiction because $d-n<d-t+1<d$ and $\fp_0\in X$.
\end{proof}
It should be mentioned that Theorem \ref{t3} is new even in the case of cyclic module. 
\begin{cor}
	Let $(R,\fm)$ be a Cohen-Macaulay local ring of dimension $d$ and let $X$ be a subset of $\Spec R$
	which is stable under generalization. Assume that $I$ and $J$ are ideals of $R$ which are linked by the zero ideal.
	If $\pd_R(R/I)<\infty$ or $R$ is Gorenstein, then the following are equivalent:
	\begin{enumerate}[(i)]
		\item{$X\subseteq S_n(R/I)$.}
		\item{$\Att_R(\hh^i_{\fm}(R/J))\subseteq\Spec R\setminus X$ for $d-n<i<d$.}
	\end{enumerate}
\end{cor}
As a consequence of Theorem \ref{t3}, we have the following result.
\begin{cor}
Let $(R,\fm)$ be a Cohen-Macaulay local ring of dimension $d$. Assume that $R_\fp$ is Gorenstein for all $\fp\in X^{n-1}(R)$.
For a horizontally linked $R$--module $M$, the following are equivalent:
\begin{enumerate}[(i)]
\item{$M_\fp$ is Cohen-Macaulay for all $\fp\in X^{n-1}(R)$.}
\item{$\depth R_\fp\geq n$ for all $\fp\in\underset{\tiny{0<i<d}}{\bigcup}\Att_R(\hh^i_{\fm}(\lambda M))$.}
\end{enumerate}
\end{cor}
\begin{proof}
The condition (i) is equivalent to say that $X^{n-1}(R)\subseteq S_d(M)$. Now the assertion is clear by Theorem \ref{t3}.
\end{proof}
The following is an immediate consequence of Theorem \ref{t3} and Theorem \ref{att}(ii).
\begin{cor}
Let $(R,\fm)$ be a Cohen-Macaulay local ring of dimension $d$. Assume that $M$ is a horizontally linked $R$--module and that $\gd_{R_\fp}(M_\fp)<\infty$ for all
$\fp\in\Spec R\setminus\{\fm\}$. The following statements are equivalent:
\begin{enumerate}[(i)]
\item{$M_\fp$ satisfies $(S_n)$ for all $\fp\in\Spec R\setminus\{\fm\}$.}
\item{$\ell(\hh^i_\fm(\lambda M)))<\infty$ for $d-n<i<d$.}
\end{enumerate}
where $\ell(-)$ denotes the length.
\end{cor}
\begin{thm}\label{c1}
	Let $(R,\fm)$ be a formally equidimensional local ring of dimension $d$ which is the homomorphic image of a Cohen-Macaulay ring. Assume that $M$ is a horizontally linked $R$--module and that $n$ is an integer with $0<n\leq d$. Let $X$ be a subset of $\Spec R$ which is stable under generalization such that $\gd_{R_\fp}(M_\fp)<\infty$ and $R_\fp$ is Cohen-Macaulay for all $\fp\in X$. 
	If $\Att_R(\hh^j_\fm(\lambda M))\subseteq\Spec R\setminus X$ for $d-n<j<d$, then $\Att_R(\hh^i_\fm(M))\subseteq\Spec R\setminus X$ for $0<i<n$.
\end{thm}
\begin{proof}
	Assume contrarily that $\fp\in\Att_R(\hh^i_\fm(M))\cap X$ for some $0<i<n$.
	By Theorem \ref{nqat}, $\fp R_\fp\in\Att_{R_\fp}(\hh^{i-\dim R/\fp}_{\fp R_\fp}(M_\fp))$.
	In particular,
	\begin{equation}\tag{\ref{c1}.1}
	\hh^{i-\dim R/\fp}_{\fp R_\fp}(M_\fp)\neq0.
	\end{equation}
	By Theorem \ref{t3}, $M_\fp$ satisfies $\widetilde{S}_n$.
	In other words, as $R_\fp$ is Cohen-Macaulay,
	\begin{equation}\tag{\ref{c1}.2}
	\depth_{R_\fp}(M_\fp)\geq\min\{\dim R_\fp, n\}.
	\end{equation}
	On the other hand, as $R$ is equidimensional and catenary,
	\begin{equation}\tag{\ref{c1}.3}
	\dim R_\fp+\dim R/\fp=d.
	\end{equation}
	If $n\geq\dim R_\fp$, then $M_\fp$ is maximal Cohen-Macaulay by (\ref{c1}.2) and so $\hh^j_{\fp R_\fp}(M_\fp)=0$ for all $j\neq\dim R_\fp$.
	Therefore, by (\ref{c1}.1), $i-\dim R/\fp=\dim M_\fp=\dim R_\fp$ and so $i=d$ by (\ref{c1}.3) which is a contradiction.
	On the other hand, if $n<\dim R_\fp$, then $\depth_{R_\fp}(M_\fp)\geq n$ by (\ref{c1}.2) and so $\hh^j_{\fp R_\fp}(M_\fp)=0$ for all $j<n$.
	In particular, $\hh^{i-\dim R/\fp}_{\fp R_\fp}(M_\fp)=0$ which is a contradiction by (\ref{c1}.1).
\end{proof}
Recall that a subset $X$ of $\Spec R$ is called \emph{specialization-closed} if every prime ideal of $R$ containing some prime ideal in $X$ belongs to $X$. The following is an immediate consequence of Theorem \ref{c1}.
\begin{cor}\label{c}
	Let $(R,\fm)$ be a formally equidimensional local ring which is the homomorphic image of a Cohen-Macaulay ring. Let $X$ be a specialization-closed subset of $\Spec R$ such that $R_\fp$ is Gorenstein for all $\fp\in\Spec R\setminus X$. Assume $M$ is a horizontally linked $R$--module. The following statements are equivalent:
	\begin{enumerate}[(i)]
		\item{$\Att_R(\hh^i_\fm( M))\subseteq X$ for $0<i<\dim R$.}
		\item{$\Att_R(\hh^i_\fm(\lambda M))\subseteq X$ for $0<i<\dim R$.}
	\end{enumerate}
\end{cor}
In the following, we extend Theorem \ref{t3} to include a change of ring.
First we need to recall some definitions and results from \cite{G}.
\begin{dfn} Let $C$ be a semidualizing $R$--module.
	An $R$--module $M$ is called \emph{$\gc$-perfect} if $\gr_R(M)=\gkd_R(M)$.  An $R$--module $M$ is called \emph{$\gc$-Gorenstein} if $M$ is $\gc$-perfect
		and $\Ext^n_R(M,C)$ is cyclic, where $n=\gkd_R(M)$. An ideal $I$ is called $\gc$-perfect (resp. $\gc$-Gorenstein) if $R/I$ is $\gc$-perfect (resp. $\gc$-Gorenstein) as $R$--module.
\end{dfn}
Note that if $I$ is a $\gc$-Gorenstein ideal of $\gc$-dimension $n$, then
$\Ext^n_R(R/I,C)\cong R/I$.

\begin{thm}\label{G2}
	\cite[Proposition 5]{G}. Let $K$ be a semidualizing $R$--module and let $I$ be a $\gkkkd$-perfect ideal. Set $C=\Ext^{\tiny{\gr(I)}}_R(R/I,K)$. Then the following statements hold.
	\begin{itemize}
		\item[(i)]{$C$ is a semidualizing $R/I$--module.}
		\item[(ii)]{If $M$ is a $R/I$--module, then $\gkkd_R(M)<\infty$ if and only if $\gkd_{R/I}(M)<\infty$, and
			if these dimensions are finite, then
			$\gkkd_R(M)=\gr(I)+\gkd_{R/I}(M)$.}
	\end{itemize}
\end{thm}
\begin{cor}
Let $(R,\fm)$ be a Cohen-Macaulay local ring and let $C$ be a
semidualizing $R$--module. Assume that $M$ is an $R$--module of
finite $\gc$-dimension which is linked by a $\gc$-Gorenstein ideal
$\fc$. Set $\overline{R}=R/\fc$. Suppose that $X$ is a subset of $\Spec \overline{R}$ which is stable under generalization.  The following are equivalent:
\begin{enumerate}[(i)]
\item{$X\subseteq\widetilde{S}_n(M)$.}
\item{$\Att_{\overline{R}}(\hh^i_{\fm}(\lambda_{\overline{R}} M))\subseteq\Spec \overline{R}\setminus X$ for $\dim\overline{R}-n<i<\dim\overline{R}$.}
\end{enumerate}
\end{cor}
\begin{proof}
As $R$ is Cohen-Macaulay and $\fc$ is $\gc$-perfect, $R/\fc$ is Cohen-Macaulay. By Theorem \ref{G2}(ii), $\gd_{R/\fc}(M)<\infty$.
Now the assertion is clear by Theorem \ref{t3}.
\end{proof}
For an $R$--module $M$, set $\cc(M)=\sup\{i<\dim_R(M)\mid\hh^i_\fm(M)\neq0\}.$
In the following, we describe the attached prime ideals of
$\hh^{\tiny{\cc(M)}}_\fm(M)$ for a non-Cohen-Macaulay horizontally
linked $R$--module $M$. The following result should be compared with \cite[Corollary 4.5]{DS1}.	
\begin{thm}\label{c2}
	Let $(R,\fm)$ be a Cohen-Macaulay local ring of dimension $d$. Assume that $M$ is an $R$--module of finite and positive Gorenstein dimension which is horizontally linked. Set $\cc=\cc(\lambda M)$.
	Then $$\Att_R(\hh_\fm^{\cc}(\lambda M))=\left\lbrace \fp\in S_{d-\cc}(M)\setminus S_{d-\cc+1}(M)\mid \depth_{R_\fp}(M_\fp)=d-\cc\right\rbrace.$$
\end{thm}
\begin{proof}
 As $\hh_\fm^{\cc}(\lambda M)$ is an Artinian $R$-module, it has a natural structure as an $\widehat{R}$-module.  More precisely, it is an Artinian $\widehat{R}$-module. By Theorem \ref{att}(iv),
	\begin{equation}\tag{\ref{c2}.1}
	\Att_R(\hh_\fm^{\cc}(\lambda M))=\left\lbrace P\cap R\mid P\in\Att_{\widehat{R}}(\hh_\fm^{\cc}(\lambda M))\right\rbrace.
	\end{equation}
	Note that by local duality theorem and Lemma \ref{l1}, $\widehat{M}$ is $(d-c)$-torsionfree $\widehat{R}$-module. It follows from Theorem \ref{att}(iii), local duality theorem and Lemma \ref{l2} that
	\begin{equation}\tag{\ref{c2}.2}
	\Att_{\widehat{R}}(\hh_\fm^{\cc}(\lambda M))=\Ass_{\widehat{R}}(\Ext^{d-\cc}_{\widehat{R}}(\widehat{(\lambda M)},\omega_{\widehat{R}}))=\Ass_{\widehat{R}}(\Ext^{d-\cc}_{\widehat{R}}(\widehat{(\lambda M)},\widehat{R})).
	\end{equation}
	On the other hand, by \cite[Proposition 9.A]{Ma1},
	\begin{equation}\tag{\ref{c2}.3}
	\Ass_{R}(\Ext^{d-\cc}_{R}((\lambda M),R))=\left\lbrace P\cap R\mid P\in\Ass_{\widehat{R}}(\Ext^{d-\cc}_{\widehat{R}}(\widehat{(\lambda M)},\widehat{R}))\right\rbrace.
	\end{equation}
It follows from (\ref{c2}.1), (\ref{c2}.2), (\ref{c2}.3) that
\begin{equation}\tag{\ref{c2}.4}
\Att_R(\hh_\fm^{\cc}(\lambda M))=\Ass_{R}(\Ext^{d-\cc}_{R}((\lambda M),R)).
\end{equation}
Set $n=d-c$. As $M$ is $n$-torsionfree, by \cite[Proposition 3.6]{AR}, there exists the following exact sequence:
\begin{equation}\tag{\ref{c2}.5}
0\longrightarrow M\longrightarrow\Omega^{n+1}\Tr\Omega^{n+1}\Tr M\oplus P\longrightarrow\Omega^{n-1}\Ext^{n+1}_R(\Tr M,R)\longrightarrow0,
\end{equation}
where $P$ is projective $R$--module. Set $X=\Omega^{n+1}\Tr\Omega^{n+1}\Tr M\oplus P$ and $Y=\Omega^{n-1}\Ext^{n+1}_R(\Tr M,R)$.
Let $\fp\in\Ass_R(\Ext^{n}_R(\lambda M,R))$. Hence $\fp\in\Ass_R(\Ext^{n+1}_R(\Tr M,R))$ and so $\fp\notin S_{n+1}(M)$ by Theorem \ref{AS}.
Note that $\depth R_\fp>n$. Therefore $\depth_{R_\fp}(Y_\fp)=n-1$ and $\depth_{R_\fp}(X_\fp)>n$. By localizing the exact sequence (\ref{c2}.5) at $\fp$, we conclude that $\depth_{R_\fp}(M_\fp)=n$.

Conversely, assume that $\fp\in S_{n}(M)\setminus S_{n+1}(M)$ and that $\depth_{R_\fp}(M_\fp)=n$. Hence by (\ref{Gdim})(ii), $\depth R_\fp>n$. By localizing the exact sequence (\ref{c2}.5) at $\fp$, we conclude that $\depth_{R_\fp}(Y_\fp)=n-1$ and so $\depth_{R_\fp}(\Ext^{n+1}_R(\Tr M,R)_\fp)=0$. In other words, $\fp\in\Ass_R(\Ext^{n+1}_R(\Tr M,R))$. Now the assertion follows from (\ref{c2}.4).
\end{proof}
For a subset $X$ of $\Spec R$, we denote by $\overline{X}$ the closure of $X$ in the Zariski topology.
\begin{cor}\label{c5}
Let $(R,\fm)$ be a Cohen-Macaulay local ring of dimension $d$. Assume that $M$ is an $R$--module of finite and positive Gorenstein dimension which is horizontally linked. For a positive integer $n<d$, we have 
$$\Att_R(\hh_\fm^{d-n}(\lambda M))\setminus\overline{\underset{\tiny{d-n<i<d}}{\bigcup}\Att_R(\hh^i_{\fm}(\lambda M))}=\left\lbrace\fp\in S_{n}(M)\setminus S_{n+1}(M)\mid \depth_{R_\fp}(M_\fp)=n\right\rbrace.$$
\end{cor}
\begin{proof}
First note that by \cite[Theorem 4.2]{DS}, for a prime ideal $\fp$ of $R$ we have  
\begin{equation}\tag{\ref{c5}.1}
\fp\in S_{n}(M)\Longleftrightarrow\hh^i_{\fp R_\fp}((\lambda M)_{\fp})=0 \text{ for } \dim R_\fp-n<i<\dim R_\fp. 
\end{equation} 
 Set $X=\left\lbrace\fp\in S_{n}(M)\setminus S_{n+1}(M)\mid \depth_{R_\fp}(M_\fp)=n\right\rbrace$ and $Y=\overline{\underset{\tiny{d-n<i<d}}{\bigcup}\Att_R(\hh^i_{\fm}(\lambda M))}$. Note that $\Spec R\setminus Y$ is stable under generalization. Hence,
by Theorem \ref{t3}, $\Spec R\setminus Y\subseteq S_{n}(M)$. Now let $\fp\in\Att_R(\hh_\fm^{d-n}(\lambda M))\setminus Y$. Therefore $\fp\in S_{n}(M)$. Also 
by Theorem \ref{nqat}, $\fp R_\fp\in\Att_{R_\fp}(\hh_{\fp R_\fp}^{\dim R_\fp-n}((\lambda M)_\fp))$. Hence, by (\ref{c5}.1), $\cc_{R_\fp}(\lambda_{R_\fp}M_\fp)=\dim R_\fp-n.$
It follows from Corollary \ref{c2} that $\fp R_\fp\notin S_{n+1}(M_\fp)$ and $\depth_{R_\fp}(M_\fp)=n$. In other words, $\fp\in X$.

Conversely, if $\fp\in X$, then $\cc_{R_\fp}(\lambda_{R_\fp}M_\fp)=\dim R_\fp-n$ by (\ref{c5}.1). It follows from Corollary \ref{c2} that $\fp R_\fp\in\Att_{R_\fp}(\hh_{\fp R_\fp}^{\dim R_\fp-n}((\lambda M)_\fp))$ and so $\fp\in\Att_R(\hh_\fm^{d-n}(\lambda M))$ by Theorem \ref{t3}. Assume contrarily that $\fp\in Y$. Hence there exists a prime ideal $\fq\subseteq\fp$ such that
$\fq\in\Att_{R}(\hh^j_{\fm}(\lambda M))$ for some integer $j$ with $d-n<j<d$. By Theorem \ref{t3}, $\fq R_\fq\in\Att_{R_\fq}(\hh_{\fq R_\fq}^{j-\dim R/\fq}((\lambda M)_\fq))$. As $j<d$, we have $j-\dim R/\fq<\dim R_\fq$. Note that $\fq\in S_{n}(M)$ and so by (\ref{c5}.1)
$j-\dim R/\fq\leq\dim R_\fq-n$. Hence $j\leq d-n$ which is a contradiction.
\end{proof}
\section{Linkage and Cohen-Macaulay modules}
The first main theorem in the theory of linkage was due to C.
Peskine and L. Szpiro. They proved that the Cohen-Macaulayness property is preserved under linkage over Gorenstein rings. Note that the above statement is no longer true if the base ring is Cohen-Macaulay but not Gorenstein. On the other hand, there is an interesting generalization of the Peskine-Szpiro Theorem, due to Martin \cite{Mar1}. Let $R$ be a Gorenstein local ring and let $I$ be Cohen-Macaulay ideal linked to $J$ by an ideal $\fa$, where $\fa$ is Cohen-Macaulay and generically Gorenstein. Then $J$ is Cohen-Macaulay if and only if $I+\fb$ is a Cohen-Macaulay ideal with the same height as $\fb$, where $\fb$ is an ideal of $R$, containing $\fa$, with $\fb/\fa$ equal to a canonical module of $\omega_{R/\fa}$ \cite[Corolary 2.11]{Mar1}.

In the theory of linkage of modules, Martsinkovsky and Strooker generalized the Peskine-Szpiro Theorem for stable modules over Gorenstein local rings \cite[Proposition 8]{MS}. In this section, we study the theory of linkage for Cohen-Macaulay modules. As a consequence, we generalize the Martin's result to the module case.

A \emph{semidualizing ideal} is an ideal that is semidualizing as an $R$-module. Let $C$ be a semidualizing $R$--module. If $R$ is generically Gorenstein, then $C$ can be identified with an ideal of $R$ \cite[Proposition 3.1]{W}. Let $R$ be a Cohen-Macaulay local ring and let $\fc$ be a  semidualizing ideal. Then $\fc$ is an ideal of height one or equals $R$ \cite[Proposition 3.2]{W}.
\begin{thm}\label{ttt}
	Let $R$ be a generically Gorenstein and let $\fc$ be a semidualizing ideal of $R$. Assume that $M$ is a horizontally linked $R$--module. Then the following statements hold:
	\begin{enumerate}[(i)]
	\item{$\gfcd_R(M)=0$ if and only if $\gfcd_R(\fc\lambda M)=0$.}	
	\item{Assume that $M$ has finite $\gfc$-dimension on $X^{n-1}(R)$. Then $M$ satisfies $\widetilde{S}_n$ if and only if $\Ext^i_R(\fc\lambda M,\fc)=0$ for $0<i<n$.}
	\item{Assume that $R$ is Cohen-Macaulay and that $\gfcd_R(M)<\infty$. Then $M$ is Cohen-Macaulay if and only if $\fc\lambda M$ is so.}
    \end{enumerate}
\end{thm}
\begin{proof}
Let $\cdots\rightarrow F_1\rightarrow F_0\overset{f}{\rightarrow} M\rightarrow0$ be the minimal free resolution of $M$. From the exact sequence $(F_0)^*\overset{g}{\rightarrow}\lambda M\rightarrow0$ we get the following exact sequence $(F_0)^*\otimes_R\fc\overset{g\otimes\fc}{\rightarrow}\lambda M\otimes_R\fc\rightarrow0$. Let $\phi=h\circ(g\otimes\fc)$, where $h$ is the natural epimorphism $\lambda M\otimes_R\fc\twoheadrightarrow\fc\lambda M$. We denote by $\psi$ the tensor evaluation isomorphism $(F_0)^*\otimes_R\fc\overset{\cong}{\rightarrow}\Hom_R(F_0,\fc)$. Dualize the free resolution of $M$, into $\fc$ to obtain the exact sequence 
\begin{equation}\tag{\ref{ttt}.1}
0\rightarrow\Hom_R(M,\fc)\rightarrow\Hom_R(F_0,\fc)\rightarrow\Hom_R(F_1,\fc)\rightarrow\Tr_{\fc}(M)\rightarrow0.
\end{equation}
Set $\lambda_{\fc}(M):=\coker(\Hom(f,\fc))$.
We denote by $\alpha$ the epimorphism
$\Hom_R(F_0,\fc)\twoheadrightarrow\lambda_{\fc}(M)$.
Now it is easy to check that $\phi\circ\psi^{-1}\circ\Hom(f,\fc)=0=\alpha\circ\psi\circ\iota$, where $\iota:\ker\phi\hookrightarrow(F_0)^*\otimes_R\fc$ is the natural inclusion.
Therefore we obtain the following commutative diagram
$$\begin{CD}
&&&&&&&&\\
\ \  &&&&0@>>>\ker\phi @>\iota>>(F_0)^*\otimes_R\fc @>\phi>>\fc\lambda M @>>>0&\\
&&&&&& @VV{\cong}V @VV{\cong}V \\
\ \ &&&&0@>>>\Hom_R(M,\fc) @>>>\Hom_R(F_0,\fc) @>>> \lambda_{\fc}(M) @>>>0&  \\
\end{CD}$$\\
with exact rows. Hence we obtain the following isomorphism:
\begin{equation}\tag{\ref{ttt}.2}
\lambda_{\fc}(M)\cong\fc\lambda M.
\end{equation} 

(i). From the exact sequence (\ref{ttt}.1) we get the following short exact sequence:
\begin{equation}\tag{\ref{ttt}.3}
0\rightarrow\lambda_{\fc}(M)\rightarrow\Hom_R(F_1,\fc)\rightarrow\Tr_{\fc}(M)\rightarrow0.
\end{equation}
As $M$ is a syzygy module, it is easy to see that
$\Ext^1_R(\Tr_{\fc}(M),\fc)=0$. Now the assertion follows from \ref{Gdim}(iv), (\ref{ttt}.2) and the exact sequence (\ref{ttt}.3).

(ii). By \cite[Proposition 2.4]{DS1}, $M$ satisfies $\widetilde{S}_n$ if and only if $\Ext^i_R(\Tr_{\fc} (M),\fc)=0$ for $1\leq i\leq n$. As $M$ is horizontally linked, it satisfies $\widetilde{S}_1$ and so
$\Ext^1_R(\Tr_{\fc}(M),\fc)=0$. The exact sequence (\ref{ttt}.3) induces the following isomorphism $\Ext^i_R(\lambda_{\fc}(M),\fc)\cong\Ext^{i+1}_R(\Tr_{\fc}M,\fc)$ for all $i>0$. Now the assertion follows from (\ref{ttt}.2).

(iii). If $M$ is Cohen-Macaulay, then $\fc\lambda M$ is Cohen-Macaulay by \ref{Gdim}(ii) and part (i).

Conversely, assume $\fc\lambda M$ is Cohen-Macaulay.
Suppose contrarily that $M$ is not Cohen-Macaulay. Hence
$\gfcd_R(M)>0$. Set $n=\min\{i>0\mid\Ext^i_R(M,\fc)\neq0\}$. As $\Ext^i_R(M,\fc)=0$ for all $0<i<n$, applying functor
$(-)^{\triangledown}=\Hom_R(-,\fc)$ to the minimal free
resolution of $M$ implies the following exact sequences:
\begin{equation}\tag{\ref{ttt}.4}
0\rightarrow\lambda_{\fc}(M)\rightarrow
(F_1)^{\triangledown}\rightarrow\cdots\rightarrow(F_n)^{\triangledown}\rightarrow\Tr_{\fc}\Omega^{n-1}(M)\rightarrow0.
\end{equation}
\begin{equation}\tag{\ref{ttt}.5}
0\longrightarrow\Ext^n_R(M,\fc)\longrightarrow\Tr_{\fc}\Omega^{n-1}(M)
\end{equation}
Now let $\fp\in\Ass_{R}(\Ext^n_R(M,\fc))$. It follows from the exact sequence (\ref{ttt}.5) that $\depth_{R_\fp}((\Tr_{\fc}\Omega^{n-1}(M)_\fp)=0$. As $\Ext^n_R(M,\fc)_\fp\neq0$, by \ref{Gdim}(ii), (iii), $\depth R_\fp-\depth (M_\fp)=\gfcdp_{R_\fp}(M_\fp)\geq n$.
Since $M$ is a first syzygy, we conclude that $\depth R_\fp>n$. It follows from the exact sequence (\ref{ttt}.4) that $\depth_{R_\fp}(\lambda_{\fc}(M)_\fp)=n$ which is a contradiction, because $\lambda_{\fc}(M)\cong\fc\lambda M$ is a maximal Cohen-Mcaulay module by our assumption.
\end{proof}
In the following, we extend Theorem \ref{ttt} to include a change of ring.
\begin{cor}\label{cc}
Let $R$ be a Cohen-Macaulay local ring and let $K$ be a semidualizing $R$--module. Assume $\fa$ is a $\gkkkd$-perfect ideal of $R$ which is generically Gorenstein. Set $C=\Ext^m_R(R/\fa, K)$, where $m=\gr(\fa)$ and write $C\cong\fc/\fa$. 
Assume that $M$ and $N$ are $R$--modules such that $M\underset{\fa}{\sim}N$.
Then the following statements hold:
\begin{enumerate}[(i)]
	\item{$M$ is $\gkkkd$-perfect if and only if $\fc N$ is $\gkkkd$-perfect.}	
	\item{Assume that $M$ has finite $\gkkkd$-dimension. 
			Then $M$ satisfies $(S_n)$ if and only if $\Ext^i_R(\fc N, K)=0$ for $m<i<n+m$. Also, $M$ is Cohen-Macaulay if and only if $\fc N$ is so.}
\end{enumerate}
\end{cor}
\begin{proof}
(i). By Theorem \ref{G2}(ii) and \cite[Lemma 3.16]{Sa}, $M$ is $\gkkkd$-perfect if and only if $\gkd_{R/\fa}(M)=0$. Note that $\gr(\fc N)=\gr(N)$. Again by using Theorem \ref{G2}(ii) and \cite[Lemma 3.16]{Sa}, we have $\fc N$ is $\gkkkd$-perfect if and only if $\gkd_{R/\fa}(\fc N)=0$. Now the assertion follows from Theorem \ref{ttt}(i).

(ii). By Theorem \ref{G2}(ii), $\gkd_{R/\fa}(M)<\infty$. As $\fa$ is $\gkkkd$-perfect, $R/\fa$ is a Cohen-Macaulay ring.
By \cite[Corollary]{G}, $\Ext^i_{R/\fa}(\fc N, C)\cong\Ext^{i+m}_R(\fc N,K)$ for all $i\geq0$.
 Now the assertion follows from Theorem \ref{ttt}(ii), (iii).
\end{proof}	
The following is a generalization of \cite[Corollary 2.11]{Mar1}.
\begin{cor}\label{ttttt}
	Let $R$ be a Cohen-Macaulay local ring and let $K$ be a semidualizing $R$--module. Assume $\fa$ is a $\gkkkd$-perfect ideal of $R$ which is generically Gorenstein. Set $C=\Ext^n_R(R/\fa, K)$, where $n=\gr(\fa)$ and write $C\cong\fc/\fa$. Assume that $M$ is an $R$--module of finite $\gkkkd$-dimension and that $N$ is a Cohen-Macaulay $R$--module which is linked to $M$ by the ideal $\fa$. Then $M$ is Cohen-Macaulay if and only if $\depth_R(N/\fc N)\geq\dim(R/\fa)-1$.
\end{cor}
\begin{proof}
	Consider the following exact sequence 
	\begin{equation}\tag{\ref{ttttt}.1}
	0\rightarrow\fc N\rightarrow N\rightarrow N/\fc N \rightarrow0.
	\end{equation}
	Now the assertion follows from Theorem \ref{ttt}(ii) and the exact sequence (\ref{ttttt}.1).	
\end{proof}
It should be noted that Theorem B is a special case of Corollary \ref{ttttt}.
\begin{cor}
	Let $R$ be a Cohen-Macaulay local ring with canonical module $\omega_R$ and let $\fa$ be a Cohen-Macaulay ideal of $R$ which is generically Gorenstein. Write $\omega_{R/\fa}\cong\fc/\fa$. Assume that $I$ and $J$ are ideals of $R$ linked by the ideal $\fa$. The following statements hold:
	\begin{enumerate}[(i)]
         \item{$R/I$ satisfies $(S_n)$ if and only if $\hh^i_{\fm}(\fc/\fc\cap J)=0$ for $\dim(R/\fa)-n<i<\dim(R/\fa)$. In particular, $R/I$ is Cohen-Macaulay if and only if $\fc/\fc\cap J$ is so.}
         \item{Assume $R/I$ is Cohen-Macaulay. Then $R/J$ is Cohen-Macaulay if and only if $\depth_R(R/(I+\fc))\geq \dim R/\fa-1$.}
    \end{enumerate}	 
\end{cor}
\begin{proof}
This follows immediately from Corollaries \ref{cc},  \ref{ttttt} and the local daulity theorem.
\end{proof}
Let $n$ be a positive integer. An $R$--module $M$ is said to be an $n$th \emph{syzygy} if there exists an exact sequence
$$0\rightarrow M\rightarrow P_{0}\rightarrow\cdots\rightarrow P_{n-1},$$ where $P_0,\cdots,P_{n-1}$ are
projective. By convention, every module is a $0$th syzygy.

\begin{thm}\label{t1}
	Let $(R,\fm)$ be a Cohen-Macaulay local ring of dimension $d$ with canonical module $\omega_R$. Assume that $M$ is a horizontally linked $R$--module and that $\gd_{R_\fp}(M_\fp)<\infty$ for all $\fp\in\Spec R-\{\fm\}$. Then the following conditions are equivalent:
	\begin{enumerate}[(i)]
		\item{$M\otimes_R\omega_R$ is maximal Cohen-Macaulay.}
		\item{$\lambda M$ is $(d+1)$th syzygy module.}
	\end{enumerate}
\end{thm}
\begin{proof}
	By \cite[Theorem 10.62]{R}, there is a third quadrant
	spectral sequence:
	$$\E^{p,q}_2=\Ext^p_R(\Tor_q^R(M,\omega_R),\omega_R)\Rightarrow\Ext^{p+q}_R(M,R).$$
	As $\gd_{R_\fp}(M_\fp)<\infty$ for all $\fp\in\Spec R-\{\fm\}$, $\Tor_i^R(M,\omega_R)$ has finite length for all $i>0$ (see \ref{ac}).
	Therefore, $\E^{p,q}_2 = 0$ when $q > 0$ and $p < d$ and so we obtain the following isomorphisms:
	\begin{equation}\tag{\ref{t1}.1}
	\Ext^i_R(M\otimes_R\omega_R,\omega_R)\cong\Ext^i_R(M,R) \text{ for all } 1\leq i\leq d.
	\end{equation}
	
	(i)$\Rightarrow$(ii). It follows from (\ref{t1}.1) and our assumption that $\Ext^i_R(M,R)=0$ for all $1\leq i\leq d$.
	As $M$ is horizontally linked, $M\cong\lambda^2M\cong\Omega\Tr\lambda M$. Therefore, $\Ext^i_R(\Tr\lambda M,R)=0$ for all $1\leq i\leq d+1$ and so
	$\lambda M$ is $(d+1)$th syzygy module by \cite[Proposition 11]{M1}.
	
	(ii)$\Rightarrow$(i). By \cite[Proposition 2.6]{DS} and our assumption, $\Ext^i_R(M,R)_\fp=0$ for all $\fp\in\Spec R-\{\fm\}$ and all $1\leq i\leq\depth R_\fp$. Hence, by (\ref{Gdim})(iii), $\gd_{R_\fp}(M_\fp)=0$ and so $\gd_{R_\fp}((\lambda M)_\fp)=0$ for all $\fp\in\Spec R-\{\fm\}$.
	It follows from \cite[Theorem 43]{M1} that $\Ext^i_R(\Tr\lambda M,R)=0$ for all $1\leq i\leq d+1$. As $M$ is horizontally linked, $M\cong\Omega\Tr\lambda M$. Therefore, $\Ext^i_R(M,R)=0$ for all $1\leq i\leq d$. Now the assertion is clear by (\ref{t1}.1).
\end{proof}
Here is a characterization of Gorenstein rings in terms of linkage.
\begin{prop}
	Let $(R,\fm,k)$ be a Cohen-Macaulay local ring of dimension $d>1$.
	Then the following statements are equivalent:
	\begin{enumerate}[(i)]
		\item{$R$ is Gorenstein.}
		\item{The Cohen-Macaulayness property is presereved under horizontal linkage.}
	\end{enumerate}
\end{prop}
\begin{proof}
(i)$\Rightarrow$(ii). Follows from \cite[Proposition 8]{MS}. (ii)$\Rightarrow$(i). By \cite[Corollary 1.2.5]{Av}, $\Omega^{d+1}k$ is stable and so it is horizontally linked (see \ref{link}). It follows from the assumption that $\lambda(\Omega^{d+1}k)$ is maximal Cohen-Macaulay. Hence $\depth_R(\Tr\Omega^{d+1}k)\geq d-1\geq1$. Now consider the following exact sequence:
$$0\rightarrow\Ext^{d+2}_R(k,R)\rightarrow\Tr\Omega^{d+1}k\rightarrow\Omega\Tr\Omega^{d+2}k\rightarrow0.$$
As $\Ext^{d+2}_R(k,R)$ has finite length and $\depth_R(\Tr\Omega^{d+1}k)>0$, we conclude that $\Ext^{d+2}_R(k,R)=0$. In other words, $R$ is Gorenstein.
\end{proof}
\section{Duality between local cohomology modules of linked module}
Recall that an $R$--module $M$ of dimension $d\geq1$ is called a
\emph{ generalized Cohen-Macaulay} module if
$\ell(\hh^i_\fm(M))<\infty$ for all $i$, $0\leq i\leq d-1$, where
$\ell$ denotes the length. 
Let $R$ be a Gorenstein local ring and let $\fa$, $\fb$ be ideals of $R$ linked by a Gorenstein ideal
$\fc$. Assume that $R/\fa$ is generalized Cohen-Macaulay. Schenzel established an interesting duality between local cohomology modules of $R/\fa$ and $R/\fb$. More precisely, he proved that $$\hh^i_{\fm}(R/\fa)\cong\Hom_R(\hh^{d-i}_{\fm}(R/\fb),\E_R(R/\fm)),$$
for  $0<i<d$, where $d=\dim R/\fa=\dim R/\fb$  \cite[Corollary 3.3]{Sc}. Martsinkovsky and Strooker generalized the above result for all generalized Cohen-Macaulay module which is linked by a Gorenstein ideal \cite[Theorem 11]{MS}. In this Section, we generalize the above result to modules which satisfy Serre's condition $(S_n)$ on the punctured spectrum of $R$.

For an ideal $\fa$ of $R$, we denote by $V(\fa)$ the set of all prime ideals of $R$ containing $\fa$. Let $X$ be a closed subset of $\Spec R$. Then $X=V(\fa)$ for some ideal $\fa$ of $R$. We call such an ideal $\fa$ the \emph{defining ideal} of $X$. This is uniquely determined up to radical.
\begin{thm}\label{t4}
Let $U$ be an open subset of $\Spec R$ and let $M$ be an $R$--module such that $\gd_{R_\fp}(M_\fp)<\infty$ for all $\fp\in U$. Assume that $U\subseteq\widetilde{S}_n(M)$ for some $0<n\leq\gr(\fa)$, where $\fa$ is the defining ideal of the complement of $U$ in $\Spec R$. Then 
$\hh^i_{\fa}(M)\cong\Ext^{i+1}_R(\Tr M,R)$ for $0\leq i<n.$
In particular, $\hh^i_{\fa}(M)$ is finitely generated for $0\leq i<n$.
\end{thm}
\begin{proof}
As $M_\fp$ satisfies $\widetilde{S}_n$,
$\Ext^i_R(\Tr M,R)_\fp=0$ for all $1\leq i\leq n$ and all
$\fp\in\Spec(R)\setminus\V(\fa)$ by Theorem \ref{AS}. Therefore,
$\underset{1\leq i\leq n}{\bigcup}\Supp_R(\Ext^i_R(\Tr M,R))\subseteq\V(\fa)$ and so $\Ext^i_R(\Tr M,R)$ is
$\fa$-torsion for $1\leq i\leq n$. Consider the following exact sequence:
\begin{equation}\tag{\ref{t4}.1}
0\longrightarrow\Ext^i_R(\Tr M,R)\longrightarrow\Tr\Omega^{i-1}\Tr M\longrightarrow
X_i\longrightarrow 0,
\end{equation}
where $X_i\approx\Omega\Tr\Omega^{i}\Tr M$. Since
$\Ext^i_R(\Tr M,R)$ is $\fa$-torsion for $1\leq i\leq n$, applying the functor $\Gamma_{\fa}(-)$ to the exact sequence (\ref{t4}.1), we get
\begin{equation}\tag{\ref{t4}.2}
\hh^j_{\fa}(\Tr\Omega^{i-2}\Tr
M)\cong\hh^j_{\fa}(\Omega\Tr\Omega^{i-1}\Tr M) \text{ for all}\  j>0,\ 2\leq i\leq n+1,
\end{equation}
and also
\begin{equation}\tag{\ref{t4}.3}
\Ext^i_R(\Tr M,R)=\Gamma_{\fa}(\Ext^i_R(\Tr M,R))\cong\Gamma_{\fa} (\Tr\Omega^{i-1}\Tr M) \text{ for all}\ 1\leq i\leq n,
\end{equation}
because $\gr(\fa)>0$.
On the other hand, if $n>1$, then
\begin{equation}\tag{\ref{t4}.4}
\hh^j_{\fa}(\Tr\Omega^{i-1}\Tr M)\cong\hh^{j+1}_{\fa}(\Omega\Tr\Omega^{i-1}\Tr M) \text{ for
all}\ i\geq 1, 0\leq j<n-1.
\end{equation}
Now by using (\ref{t4}.2), (\ref{t4}.3) and (\ref{t4}.4) we obtain
the result.
\end{proof}
The following is an immediate consequence of Theorem \ref{t4}.
\begin{cor}\label{c4}
Let $U$ be an open subset of $\Spec R$ and let $M$ be a horizontally linked $R$--module such that $\gd_{R_\fp}(M_\fp)<\infty$ for all $\fp\in U$. Assume that $U\subseteq\widetilde{S}_n(M)$ for some $1<n\leq\gr(\fa)$, where $\fa$ is the defining ideal of the complement of $U$ in $\Spec R$. Then
$\hh^i_{\fa}(M)\cong\Ext^i_R(\lambda M,R) \text{ for }  0<i<n.$
\end{cor}
Let $R$ be a Cohen-Macaulay local ring. An $R$--module $M$ is generalized Cohen-Macaulay if and only if $M_\fp$ is a Cohen-Macaulay $R_\fp$--module for all $\fp\in\Spec(R)-\{\fm\}$ (see \cite[Lemma 1.2 , Lemma 1.4]{T}). Therefore the following result can be viewed as a generalization of \cite[Theorem 10]{MS}.
\begin{cor}\label{c6}
Let $(R,\fm,k)$ be a Gorenstein local ring of dimension $d>1$ and let
$M$ be a horizontally linked $R$--module. Assume $n$ is a positive
integer and that $M_\fp$ satisfies $(S_n)$ for all
$\fp\in\Spec R\setminus\{\fm\}$. Then
$\hh^i_{\fm}(M)\cong\Hom_R(\hh^{d-i}_{\fm}(\lambda M),\E_R(k))
\text{ for } 0<i<n$, where $\E_R(k)$ is the injective envelope of $k$.
\end{cor}
\begin{proof}
Follows from Corollary \ref{c4} and the local duality theorem.	
\end{proof}
We end this section by the following result which extend Corollary \ref{c6} to include a change of ring.
\begin{cor}
Let $(R,\fm,k)$ be a Gorenstein local ring, $\fc$ a Gorenstein ideal of $R$, $\overline{R}=R/\fc$, and $d=\dim\overline{R}>1$.
Assume that $M$ is an $R$--module which is linked by the ideal $\fc$ and that $M_\fp$ satisfies $(S_n)$ for all
$\fp\in\Spec \overline{R}\setminus\{\overline{\fm}\}$. Then
$\hh^i_{\fm}(M)\cong\Hom_R(\hh^{d-i}_{\fm}(\lambda_{\overline{R}} M),\E_R(k))
\text{ for } 0<i<n$.
\end{cor}
\bibliographystyle{amsplain}

\begin{thebibliography}{9}

\bibitem{AF}
~F.~W. Anderson and ~K. ~R. Fuller, \emph{Rings and Categories of
Modules}, Second edition, Springer-Verlag, 1992.

\bibitem{AB}
~M. Auslander and ~M. Bridger, \emph{Stable module theory}, Mem. of
the AMS 94, Amer. Math. Soc., Providence 1969.

\bibitem{AR}
~M. Auslander and ~I. Reiten.
\newblock Syzygy modules for Noetherian rings.
\newblock J. Algebra {\bf 183} (1996), no. 1, 167--185.

\bibitem{Av}
~L.~L. Avramov, \emph{Infinite free resolutions}, Six lectures on commutative algebra (Bellaterra, 1996), Progress in Mathematics, vol. 166,
Birkh$\ddot{a}$user, Basel, 1998, pp. 1--118.

\bibitem{BS}
~M.~P. Brodmann and R.~Y. Sharp, \emph{Local cohomology: an
algebraic introduction with geometric applications}, Cambridge
Studies in Advanced Mathematics, 60. Cambridge University Press,
Cambridge, 1998.

\bibitem{CDGST}
~O. Celikbas, M. T. Dibaei, M. Gheibi, A. Sadeghi and R. Takahashi, \emph{Associated Primes and Syzygies of Linked Modules}. J. Commut. Algebra (To appear).

\bibitem{DGHS}
~M. T. Dibaei, ~M. Gheibi, ~S. ~H. Hassanzadeh and ~A. Sadeghi, \emph{Linkage of modules over Cohen-Macaulay rings}, J. Algebra, 335 (2011) 177--187.

\bibitem{DS}
~M.~T. Dibaei,~A. Sadeghi, \emph{Linkage of finite Gorenstein dimension modules }, J. Algebra 376 (2013) 261--278.

\bibitem{DS1}
~M.~T. Dibaei,~A. Sadeghi, \emph{Linkage of modules and the Serre conditions}, J. Pure Appl. Algebra 219 (2015) 4458--4478.

\bibitem{DS2}
~M.~T. Dibaei,~A. Sadeghi, \emph{Linkage of modules with respect to a semidualizing module}, Pacific J. Math. 294 (2018), No. 2, 307--328.

\bibitem{F}
~H.~B. Foxby, \emph{Gorenstein modules and related modules}, Math. Scand. 31 (1972), 267--285.

\bibitem{F1}
~H.~B. Foxby, \emph{Quasi-perfect modules over Cohen-Macaulay
rings}, Math. Nachr. 66 (1975), 103--110.

\bibitem{G}
~E.~S. Golod, \emph{G-dimension and generalized perfect ideals},
Trudy Mat. Inst. Steklov. 165 (1984), 62-–66; English transl. in
Proc. Steklov Inst. Math. 165 (1985).

\bibitem{HVV} ~J. Herzog, ~W. ~V. Vasconcelos and ~R. Villarreal, \emph{ Ideals with sliding depth}, Nagoya Math.
J. 99 (1985), 159--172.

\bibitem{Hu}
~C. Huneke, \emph{Linkage and the Koszul homology of ideals}, Amer. J. Math. 104 (1982), no. 5, 1043--1062.

\bibitem{IT}
~K.-i. Iima and ~R. Takahashi, \emph{Perfect linkages of modules}, J. Algebra 458 (2016), 134--155.

\bibitem{IGM}
~I.~G. Macdonald, \emph{Secondary representation of modules over a commutative ring}, Symp. Math. 11 (1973) 23--43.

\bibitem{Mar1}
~H.~M. Martin, \emph{Linkage by generically Gorenstein, Cohen-Macaulay ideals.}
J. Algebra 207 (1998), no. 1, 43--71.

\bibitem{Mar}
~H.~M. Martin, \emph{Linkage and the generic homology of modules}, Comm. Algebra 28 (2000) 199--213.

\bibitem{MS}
~A. Martsinkovsky and ~J.~R. Strooker, \emph{Linkage of modules}, J. Algebra 271 (2004), 587--626.

\bibitem{M1}
~V. Ma\c{s}iek, \emph{Gorenstein dimension and torsion of modules over commutative Noetherian rings}, Comm. Algebra (2000), 5783--5812.

\bibitem{Ma1}
~H. Matsumura, \emph{Commutative algebra}, Second Edition, Benjamin 1980.

\bibitem{Ma}
~H. Matsumura, \emph{Commutative ring theory}, Cambridge University
Press 1986.

\bibitem{Na}
~U. Nagel, \emph{Liaison classes of modules}, J. Algebra 284 (2005), no. 1, 236--272.

\bibitem{NQ}
~L.T. Nhan and ~P.H. Quy, \emph{Attached primes of local cohomology modules under localization and completion}, J.
Algebra 420 (2014), 475--485.

\bibitem{N}
~K. Nishida, \emph{Linkage and duality of modules}, Math. J. Okayama Univ. 51 (2009), 71--81.

\bibitem{PS}
~C. Peskine and ~L. Szpiro, \emph{Liasion des vari\'{e}t\'{e}s
alg\'{e}briques.} I, Inv. math, 26 (1974), 271--302.

\bibitem{TJP}
~T.~J. Puthenpurakal, \emph{A function on the set of isomorphism classes in the stable category of maximal Cohen-Macaulay modules over a Gorenstein ring: with applications to liaison theory}, Math. Scand. 120 (2017), 161--180.

\bibitem{TP}
~T.~J. Puthenpurakal, \emph{Invariants of linkage of modules}, Preprint 2015; posted at arXiv:1512.05105.

\bibitem{R}
~J. Rotman, \emph{An Introduction to Homological Algebra}, Academic
Press, New York (1979).

\bibitem{Sa}
~A. Sadeghi, \emph{Linkage of finite $G_C$-dimension modules}. J. Pure Appl. Algebra 221 (2017) 1344--1365.

\bibitem{Sc}
~P. Schenzel, \emph{Notes on liaison and duality }, J. Math. Kyoto
Univ. 22 (1982/83), no. 3, 485--498.

\bibitem{T}
~N. ~V. Trung, \emph{Toward a theory of generalized Cohen-Macaulay modules}, Nagoya Math. J.
Vol. 102 (1986), 1--49.

\bibitem{W}
S. Sather-Wagstaff, \emph{Semidualizing modules and the divisor class group}, Illinois J. Math. 51-1 (2007), 255--285.

\bibitem{YI}
~Y. Yoshino and ~S. Isogawa, \emph{Linkage of Cohen-Macaulay modules over a Gorenstein ring}, J. Pure Appl. Algebra 149
(2000), 305--318.

\end{thebibliography}

\end{document}